%Cgeode.tex: 
%%a Plain TeX file by Tewodros Amdeberhan, Manuel Kauers, and  Doron Zeilberger (x pages)

%begin macros

\baselineskip=14pt
\parskip=10pt

\font\eightrm=cmr8 

\magnification=\magstephalf

\def\1{{\overline{1}}}
\def\2{{\overline{2}}}
\parindent=0pt
\overfullrule=0in

\def\frac#1#2{{#1 \over #2}}
%\headline={\rm  \ifodd\pageno  \RightHead  \else  \LeftHead  \fi}
%\def\RightHead{\centerline{
%Title
%}}
%\def\LeftHead{ \centerline{Doron Zeilberger}}
%end macros
\centerline
{\bf  The Challenge of Computing Geode Numbers }

\bigskip
\centerline
{\it Tewodros AMDEBERHAN, Manuel KAUERS, and Doron ZEILBERGER}

{\bf Abstract}:
In a fascinating recent American Mathematical Monthly article, Norman Wildberger and Dean Rubine introduced a new kind of combinatorial numbers, that they aptly named the ``Geode numbers''.
While their definition is simple, these numbers are surprisingly hard to compute, in general.
While the two-dimensional case has a nice closed-form expression, that make them easy to compute,
already the three-dimensional case poses major computational challenges that we do meet, combining experimental mathematics and the holonomic ansatz.
Alas, things get really complicated in four and higher dimensions, and
we are unable to efficiently compute, for example, the $1000$-th term of the four-dimensional diagonal Geode sequence. A donation of $100$ US dollars to the OEIS, in honor of the first person to compute this number, is offered.
\bigskip

{\bf Maple package:} This article is accompanied by the Maple package {\tt Geode.txt}, available directly from the url : \quad{\tt https://sites.math.rutgers.edu/\~{}zeilberg/tokhniot/Geode.txt} \quad.

There are numerous input and output files in the front of this article:

{\tt https://sites.math.rutgers.edu/\~{}zeilberg/mamarim/mamarimhtml/Cgeode.html} \quad .

{\bf Introduction}

In the May 2025 issue of the {\it American Mathematical Monthly}, Dean Rubine and Norman Wildberger, in the course of stating their beautiful
``explicit'' solution to the general algebraic equation (with {\it symbolic} coefficients $c_0,c_1,c_2, \dots, c_k$)
$$
0= c_0 -c_1 \alpha +c_2 \alpha^2+ \dots + c_k \alpha^k \quad ,
$$
needed the following numbers defined on lists  of non-negative integers $[m_1,m_2, \dots, m_k ]$,
that they called the {\it hyper-Catalan numbers}:
$$
C(m_1, \dots, m_k):=\frac{(2m_1+ 3m_2+ \dots + (k+1)m_k)!}{(1+m_1+ 2m_2+ \dots + km_k)!\, m_1! \cdots m_k!} \quad .
$$

Let $n$ and $k$ be positive integers. Define a polynomial, let's call it $P_{n,k}(t_1, \dots , t_k)$, in the $k$ (continuous) variables, $t_1, \dots, t_k$, of total degree $n$, as follows:
$$
P_{n,k} (t_1, \dots, t_k):=\sum_{m_1+m_2+ \dots + m_k=n} \, C(m_1, \dots, m_k) t_1^{m_1} \cdots t_k^{m_k} \quad.
$$

Rubine and Wildberger proved the following surprising theorem (Theorem 12 of [WR]).

{\bf The Geode Theorem} (Rubine and Wildberger): For all positive integers $n$ and $k$, the polynomial
$P_{n,k}(t_1, \dots, t_k)$ is divisible by $t_1+\dots+ t_k$.

So it makes sense to define another polynomial, of degree $n$ in $t_1, \dots, t_k$:

{\bf Definition} (The Geode polynomial). For positive integers $n$ and $k$
$$
Q_{n,k}(t_1, \dots, t_k):= \frac{P_{n+1,k}(t_1,\dots, t_k)}{t_1+ \dots + t_k} \quad \quad .
$$

The {\it geode} numbers are the coefficients of this polynomial.

{\bf Definition} (The Geode numbers) for every list of non-negative integers, $[m_1, \dots, m_k]$,
the coefficient of $t_1^{m_1} \cdots t_k^{m_k}$ in $Q_{m_1+ \dots + m_k,k}(t_1, \dots, t_k)$ is $G(m_1, \dots, m_k)$.

{\eightrm Procedure {\tt G(m)} in our Maple package uses the above definition to compute these numbers.}

For example, in order to find $G(4,7,8)$ type {\tt G([4,7,8]);}, getting $11258614474275030033600$. This works for small arguments, but
don't even try to do {\tt G([1000,1000,1000]);} . It will take your computer for ever. Since the computer would first
compute the (total) degree-$3003$ polynomial  $P_{3003,3}(t_1,t_2,t_3)$, then divide by $t_1+t_2+t_3$, and then extract the coefficient
of $t_1^{1000}\,t_2^{1000}\,t_3^{1000}$. As we will soon see, and that is one of the achievements of the present paper, using
experimental mathematics and {\it guessing}, this $3910$-digit integer can be computed in $0.231$ seconds. Alas, at this time of writing
we can't compute $G(1000,1000,1000,1000)$ and one of us (DZ) is pledging to donate \$100 dollars to  {\tt http://oeisf.org/}
in honor of the first person to compute it correctly.

{\bf The Two-Dimensional Geode Numbers}

It turned out that $G(m_1,m_2)$ has a nice {\bf closed-form} expression, very similar to $C(m_1,m_2)$.
It was conjectured in [WR] and soon proved (along with some other conjectures from that paper) in [AZ] and [R].

{\bf Theorem 1} (Conjectured in [WR] and proved in [AZ] and [R])
$$
G(m_1,m_2) \, = \, \frac{1}{(2m_1+2m_2+3)(m_2+m_3+1)} \cdot \frac{(2m_1+3m_2+3)!}{(m_1+2m_2+2)!m_1!m_2!} \quad.
$$

{\eightrm This is implemented in procedure {\tt G2f(m)} in our Maple package, that inputs $m=[m_1,m_2]$ and
outputs $G(m_1,m_2)$. For example, computing $G(5000,5000)$ using the definition (i.e. executing {\tt G([5000,5000]);})
takes on our laptop $91$ seconds, while using the above explicit formula only takes $0.005$ seconds!}

{\bf The Three-Dimensional Geode Numbers}

An equivalent way of stating theorem 1, is saying that the discrete bi-variate function, $G(m_1,m_2)$ satisfies the
following {\it pure}, {\bf first-order} linear recurrences
$$
G(m_1,m_2)=
\frac{\left(2 {m_1} +2+3 {m_2} \right) \left(2 {m_1} +3{m_2} +3\right) \left(2 {m_1} +1+2{m_2} \right) \left({m_1} +{m_2} \right)}{{m_1} \left(2 {m_1} +2 {m_2} +3\right)
\left({m_1} +{m_2} +1\right) \left({m_1} +2 {m_2} +2\right)} \cdot G(m_1-1,m_2) \quad ,
$$
\vfill\eject
$$
G(m_1,m_2)=
$$
$$
\frac{\left(2 {m_1} +1+3{m_2} \right) \left(2 {m_1} +2+3{m_2} \right) \left(2 {m_1} +3{m_2} +3\right) \left(2 {m_1} +1+2{m_2} \right) \left({m_1} +{m_2} \right)}
{{m_2} \left(2 {m_1} +2 {m_2} +3\right) \left({m_1} +{m_2} +1\right) \left(1+{m_1} +2 {m_2} \right) \left({m_1} +2{m_2} +2\right)}
\cdot G(m_1,m_2-1) \quad .
$$

Together with the {\it initial condition} $G(1,1)=16$, this enables one to compute any $G(m_1,m_2)$ very fast, in fact it is more efficient to use
the equivalent recurrence than the ``closed-form'' expression for $G(m_1,m_2)$ given by Theorem $1$, since you divide very large integers by other very
large integers.

While there is no `closed-form' expression for $G(m_1,m_2,m_3)$, the next-best-thing is true!
$G(m_1,m_2,m_3)$ satisfies {\bf second-order}, pure, linear recurrences with polynomial coefficients, in each of its arguments.
This also enables very fast
computation of $G(m_1,m_2,m_3)$ for large arguments, that would be impractical using the definition.

So we have the next very useful theorem.

{\bf Theorem 2}: There exist (explicit) rational functions in $m_1,m_2,m_3$,

$\bullet$ $f_{1,1}(m_1,m_2,m_3)$ and $f_{1,2}(m_1,m_2,m_3)$ both with numerators and denominators of (total) degree $11$

$\bullet$ $f_{2,1}(m_1,m_2,m_3)$ and $f_{2,2}(m_1,m_2,m_3)$ both with numerators and denominators of (total) degree $14$

$\bullet$ $f_{3,1}(m_1,m_2,m_3)$ and $f_{3,2}(m_1,m_2,m_3)$ both with numerators and denominators of (total) degree $17$

such that $G(m_1,m_2,m_3)$ satisfies the second-order linear recurrences
$$
G(m_1,m_2,m_3) \, = \, f_{1,1}(m_1,m_2,m_3) G(m_1-1,m_2,m_3) \, + \,f_{1,2}(m_1,m_2,m_3) G(m_1-2,m_2,m_3)  \quad ,
$$
$$
G(m_1,m_2,m_3) \, = \, f_{2,1}(m_1,m_2,m_3) G(m_1,m_2-1,m_3) \, + \,f_{2,2}(m_1,m_2,m_3) G(m_1,m_2-2,m_3)  \quad ,
$$
$$
G(m_1,m_2,m_3) \, = \, f_{3,1}(m_1,m_2,m_3) G(m_1,m_2,m_3-1) \, + \,f_{3,2}(m_1,m_2,m_3) G(m_1,m_2,m_3-2)  \quad .
$$

{\eightrm These rational functions are too complicated to be presented here, but are easily viewable via our Maple package by
typing

{\tt Key3(m1,m2,m3)[1][1], Key3(m1,m2,m3)[1][2]},  {\tt Key3(m1,m2,m3)[2][1], Key3(m1,m2,m3)[2][2]},
{\tt Key3(m1,m2,m3)[3][1], Key3(m1,m2,m3)[3][2]},  respectively. }

Note that together with the initial conditions $G(i,j,k)$ with $1 \leq i,j,k \leq 2$, these enable
very fast computation of $G(m_1,m_2,m_3)$.

If you are only interested in the diagonal sequence $G(n,n,n)$, then we have the following deep theorem.

{\bf Theorem 3}: There exist rational functions in $n$, $f_1(n)$ and $f_2(n)$,  whose numerators and denominators
are polynomials of degree $35$ in $n$, such that the sequence $\{G(n,n,n)\}_{n=1}^{\infty}$ satisfies the
second-order linear recurrence
$$
G(n,n,n) \, = \, f_1(n) \, G(n-1,n-1,n-1) \,+ \, f_2(n) \, G(n-2,n-2,n-2) \quad.
$$

Again $f_1(n)$ and $f_2(n)$ are too complicated to be displayed here, but you can see them by typing

{\tt Rec3AllPC(n)[1][1][1]}, {\tt Rec3AllPC(n)[1][1][2]},  respectively.

Note that together with the initial conditions $G(1,1,1)=319$ and $G(2,2,2)=669123$ this recurrence allows us to compute, very fast,
the values of $G(n,n,n)$ for large $n$ that would he hopeless using the definition.

{\eightrm This is implemented in procedure {\tt G3diag(n)}. To see the values of
$G(10^i,10^i,10^i)$ for $1 \leq i \leq 5$, see the output file:

{\tt https://sites.math.rutgers.edu/\~{}zeilberg/tokhniot/oGeode2.txt} }.

{\bf How did we find these amazing (and useful!) recurrences for the 3D Geode numbers?}

By guessing, of course!

We cranked out sufficiently many terms using the definition
of $G(m_1,m_2,m_3)$ and then used {\it linear algebra} and {\it undetermined coefficients}
to guess the recurrences, just like in {\it machine learning}. We had the advantage that
the data was {\it exact} without noise. Once the computer found the recurrenecs using
the {\it training data}, we could verify it for many other cases.

Once guessed it is routine (and straightforward) to prove our guesses fully rigorously, as follows.

By using the general theory [Z] [Ka], it is readily seen that $G(m_1,m_2,m_3)$, from its definition, is
holonomic, i.e. there exist linear recurrences with polynomial coefficients, in each of the directions, that it satisfies.
The class of holonomic functions is an {\it algebra}. Let $G'(m_1,m_2,m_3)$ be the three-variate
discrete function {\it defined} by the above recurrences (of course one would want to verify that the recurrences are
compatible, but  this is easy). We have to prove that for all $m_1, m_2,m_3 \geq 0$,
$$
G(m_1,m_2,m_3)=G'(m_1,m_2,m_3) \quad.
$$

An equivalent way of defining $G(m_1,m_2,m_3)$ is as the unique function on ${\bf N}^3$ satisfying
$$
G(m_1-1,m_2,m_3)+G(m_1,m_2-1,m_3)+G(m_1,m_2,m_3-1) \, = \, C(m_1,m_2,m_3) \quad.
$$
Let $F(m_1,m_2,m_3)$ be defined by
$$
F(m_1,m_2,m_3):=G'(m_1-1,m_2,m_3)+G'(m_1,m_2-1,m_3)+G'(m_1,m_2,m_3-1)- C(m_1,m_2,m_3) \quad.
$$
In order to prove that $G(m_1,m_2,m_3)=G'(m_1,m_2,m_3)$,  for all $m_1,m_2,m_3 \geq 0$, we need to prove that
$$
F(m_1,m_2,m_3) \, = \, 0 \quad,
$$
for all $m_1,m_2,m_3 \geq 0$.
Being a linear combination of holonomic functions $F(m_1,m_2,m_3)$ is automatically also holonomic, and
there are standard algorithms, implemented in Christoph Koutschan's {\it holonomic calculator} Mathematica package
that do it. But why bother? Since we know beforehand that $F(m_1,m_2,m_3)=0$ for all $m_1 \leq K, m_2 \leq K, m_3 \leq K$ for
some easily determined $K$, and after checking that the orders are $\leq K$, it follows immediately that $F(m_1,m_2,m_3)$ is
identically zero. Giving a {\it rigorous} proof to those obtuse readers who insist on it.

{\bf How about 4 Dimensions and beyond?}

We were unable to guess a linear recurrence with polynomial coefficients for $\{G(n,n,n,n)\}$, or pure recurrences for $G(m_1,m_2,m_3,m_4)$. While
we know that they {\bf exist} for sure, finding them is a different matter.

One of us (DZ) is offering to donate $100$ US dollars to the OEIS for the determination of
$$
G(1000,1000,1000,1000) \quad,
$$

and $200$ US dollars for
$$
G(1000,1000,1000,1000,1000) \quad.
$$

Good luck!

{\bf References}

[AZ] Tewodros Amdeberhan and Doron Zeilbeger, {\it Proofs of Three Geode Conjectures}, submitted. \hfill\break
{\tt https://sites.math.rutgers.edu/\~{}zeilberg/mamarim/mamarimhtml/geode.html} \hfill\break
{\tt https://arxiv.org/abs/2506.17862} \quad.

[Ge] Ira Gessel, {\it Lattice paths and the Geode}, {\tt https://arxiv.org/abs/2507.09405} \quad .

[Go] Fern Gossow, {\it Ordered trees and the Geode}, {\tt https://arxiv.org/abs/2507.18097} \quad .

[Ka] Manuel Kauers, {\it `` D-finite functions''}, Springer, 2023.

[Ko] Christoph Koutschan. {\it Holonomic Functions (a Mathematica package)}. \hfill\break
{\tt https://risc.jku.at/sw/holonomicfunctions/}.

[R] Dean Rubine, {\it Hyper-Catalan and Geode Recurrences and Three Conjectures of Wildberger}, 
{\tt https://arxiv.org/abs/2507.04552} \quad.

[WR] N. J. Wildberger and  D.Rubine, {{\it Hyper-Catalan Series Solution to Polynomial Equations, and the Geode},
 Amer. Math. Monthly, {\bf 132} (2025),  383-402. \hfill\break
{\tt https://www.tandfonline.com/doi/full/10.1080/00029890.2025.2460966}

[Z] Doron Zeilberger, {\it A Holonomic Systems Approach To Special Functions},
 J. Computational and Applied Math {\bf 32} (1990), 321-368. \hfill\break
{\tt https://sites.math.rutgers.edu/\~{}zeilberg/mamarim/mamarimhtml/holonomic.html} \quad .

\bigskip
\hrule
\bigskip
Tewodros Amdeberhan, Department of Mathematics, Tulane University, New Orleans, LA 70118, USA; 
{\tt tamdeber@tulane.edu} \quad .
\smallskip
Manuel Kauers, Institute for Algebra, J. Kepler University Linz, Austria; \hfill\break
{\tt manuel.kauers@jku.at} $\,\,$ .
\smallskip
Doron Zeilberger, Department of Mathematics, Rutgers University (New Brunswick), Hill Center-Busch Campus, 110 Frelinghuysen
Rd., Piscataway, NJ 08854, USA;
{\tt DoronZeil@gmail.com} $\,\,$.

\bigskip

Aug. 13, 2025.

\end